\documentclass[11pt]{article}

\usepackage{fullpage}
\usepackage{math}
\usepackage[colorlinks=true,linkcolor=blue,urlcolor=black,citecolor=blue,breaklinks]{hyperref}
\usepackage[T1]{fontenc}
\usepackage{color}
\usepackage{enumerate}
\usepackage{anyfontsize}
\usepackage{cite}
\usepackage[ruled,norelsize]{algorithm2e}
\usepackage{graphicx}
\usepackage[caption=false,font=footnotesize]{subfig}

\numberwithin{equation}{section}

\newcommand*\samethanks[1][\value{footnote}]{\footnotemark[#1]}

\begin{document}

\title{Minimax Lower Bounds for $\hinf$-Norm Estimation}

\author{Stephen Tu\thanks{Both authors contributed equally to this work.}~, Ross Boczar\samethanks[1]~, and Benjamin Recht
\\
University of California, Berkeley%
}
\maketitle

\begin{abstract}
The problem of estimating the $\hinf$-norm of an LTI system from
noisy input/output measurements has attracted recent attention
as an alternative to parameter identification
for bounding unmodeled dynamics in robust control.
In this paper, we study lower bounds for $\hinf$-norm estimation
under a query model where at each iteration the algorithm chooses
a bounded input signal and receives the response of the chosen signal
corrupted by white noise.
We prove that when the underlying system is an FIR filter,
$\hinf$-norm estimation is no more efficient than model identification
for \emph{passive} sampling.
For \emph{active} sampling, we show that norm estimation
is at most a factor of $\log{r}$ more sample efficient
than model identification, where $r$ is the length of
the filter.
We complement our theoretical results with experiments which
demonstrate that a simple non-adaptive estimator of the norm
is competitive with state-of-the-art adaptive norm estimation algorithms.
\end{abstract}

\section{Introduction}

Recently, many researchers have proposed algorithms
for estimating the $\hinf$-norm of a linear time-invariant (LTI) filter from input/output
data \cite{vanheusden07,wahlberg2010non,rojas2012analyzing,oomen14,muller17bandit,rallo17datadriven}. A common property of these algorithms
is eschewing model parameter estimation for directly estimating
either the worst case $\ell_2$-input signal \cite{wahlberg2010non,rojas2012analyzing} or the maximizing frequency \cite{muller17bandit,rallo17datadriven}.
One of the major motivations behind these algorithms is sample efficiency:
since the $\hinf$-norm is a scalar estimate whereas the number of model parameters
can be very large and possibly infinite, intuitively one expects
that norm estimation can be performed using substantially fewer samples compared with model estimation.

In this paper, we study the fundamental limits of estimating the $\hinf$-norm
of a finite impulse response (FIR) filter, and compare to known bounds for FIR
model estimation. We show that for \emph{passive} algorithms which do not adapt
their inputs in response to the result of previous queries, it is no more
efficient to estimate the $\hinf$-norm than to estimate the model. For
\emph{active} algorithms which do adapt their inputs, we show that compared to
model estimation, norm estimation is only at most a $\log{r}$ factor more
efficient when the underlying model is a length $r$ FIR filter.  Our analysis
raises an interesting open question: whether or not there there exists an
active sampling strategy which attains our stated lower bound.

Based on our theoretical findings, we study the empirical performance of
several existing adaptive algorithms compared to a simple (non-adaptive) 
estimator which first fits a model via least-squares and then returns the norm of the model.  
Surprisingly, we find that the current adaptive methods
do not perform significantly better than the simple estimator. 

\section{Related Work}

Data-driven methods for estimating the $\hinf$-norm of an LTI system
fall into roughly two major approaches: (a) estimating the worst case $\ell_2$-signal
via a power-iteration type algorithm \cite{wahlberg2010non,rojas2012analyzing,oomen14} and (b) discretizing the interval $[0, 2\pi)$ and
searching for the maximizing frequency \cite{muller17bandit,rallo17datadriven}.

Algorithms that rely on power iteration take advantage of a clever
time-reversal trick introduced by Wahlberg et al.~\cite{wahlberg2010non}, which
allows one to query the adjoint system $G^*$ with only input/output access to
the original system $G$.  One issue with these methods is that the rate
of convergence of the top singular value of the truncated Toeplitz matrix to
the $\hinf$-norm of the system is typically $O(1/n^2)$ (c.f. \cite{bottcher2000book}), but the constant
hidden in the $O(\cdot)$ notation can be quite large as pointed out by
\cite{tu18toeplitz}.  A non-asymptotic analysis of the statistical quality of
the norm estimate returned by the power iteration remains an open question;
asymptotic results can be found in \cite{rojas2012analyzing}.

The algorithms developed in \cite{muller17bandit,rallo17datadriven} 
based on discretizing the frequencies $[0, 2\pi)$
are rooted in ideas from the multi-armed bandit literature. Here, each frequency
can be treated as either an ``arm'' or an ``expert'', and an adaptive algorithm
such as Thompson sampling~\cite{muller17bandit} or multiplicative weights~\cite{rallo17datadriven}
is applied. While sharp regret analysis for bandits has been developed
by the machine learning and statistics communities~\cite{bubeck12}, one of the barriers to applying
this analysis is a lack of a sharp theory for the level of discretization required.
In practice, the number of grid points is a parameter that must be appropriately tuned
for the problem at hand.

The problem of estimating the model parameters of a LTI system with model error
measured in $\hinf$-norm is studied by \cite{tuboczar}, and in $\ell_p$-norms
by \cite{goldenshluger98}. For the FIR setting, \cite{tuboczar} gives matching
upper and lower bounds for identification in $\hinf$-norm.  These bounds will
serve as a baseline for us to compare our bounds with in the norm estimation
setting.
Helmicki et al.~\cite{helmicki91} provide lower bounds for estimating
both a model in $\hinf$-norm and its frequency response at a particular frequency
in a query setting where the noise is worst-case. In this work we consider 
a less conservative setting with stochastic noise.
M{\"{u}}ller et al.~\cite{muller17bandit} prove an asymptotic regret lower bound
over algorithms that sample only one frequency at every iteration.
Their notion of regret is however defined with respect to the best frequency
in a fixed discrete grid and not the $\hinf$-norm.
As we discuss in Section~\ref{sec:active_lb_discussion}, this turns out to be a subtle but important distinction.

\section{Problem Setup and Main Results}

In this section we formulate the problem under consideration and state our main results.
We fix a filter length $r$ and consider an unknown length $r$ causal FIR filter $H(g) := \sum_{k=0}^{r-1} g_k z^{-k}$ with $g \in \C^r$.
We study the following time-domain input/output query model for $H(g)$: for $N$ rounds, we
first choose an input $u_t \in \C^r$ such that $\norm{u_t}_2 \leq M$, and then
we observe a sample $y_t \sim \calN(T(g) u_t, \sigma^2 I)$,
where $T(g)$ denotes the $r \times r$ upper left section of the
semi-infinite Toeplitz matrix induced by treating $g$ as an element of $\ell_2$
\footnote{
We note that our results extend naturally to the setting when $T(g)$ is the $\alpha r \times \alpha r$ upper left section
for a positive integer $\alpha \geq 1$.
Furthermore, one can restrict both the system coefficients $g$ and the inputs $u_t$ to be real-valued
by considering the discrete cosine transform (DCT) instead of the discrete Fourier transform (DFT)
in our proofs.}.
By $\calN(\mu, \sigma^2 I)$ for a complex $\mu \in \C^r$ we mean we observe
$(\Re(\mu) + \xi_1) + j(\Im(\mu) + \xi_2)$ where $\xi_i \sim \calN(0, \sigma^2 I)$ and are independent.

After these $N$ rounds,
we are to return an estimate $\widehat{H} \in \R$
of the norm $\hinfnorm{H(g)}$
based on the collected data $(u_1, y_1, ..., u_N, y_N)$.
The expected risk of any algorithm for this problem is measured as $\E[\abs{\widehat{H} - \hinfnorm{H(g)}}]$, where
the expectation is taken with respect to both the randomness of the algorithm and the noise of the outputs $y_t$.

Our results distinguish between \emph{passive} and \emph{active} algorithms. A passive algorithm is one
in which the distribution of the input $u_t$ at time $t$ is
independent of the history $(u_1, y_1, ..., u_{t-1}, y_{t-1})$.
An active algorithm is one where the distribution of $u_t$ is allowed to depend
on this history.

Given this setup, our first result is a minimax lower bound for the risk attained by any passive algorithm.
\begin{mythm}[Passive lower bound]
\label{thm:passive_lower_bound}
Fix a $\gamma > 0$.
Let $r \geq c$ for a universal constant $c >0$ and $N \geq \mathrm{poly}(r, M, 1/\gamma)$. We call a passive algorithm $\calA$ admissible
if the matrix $\frac{1}{N} \sum_{t=1}^{N} \E_{u_t}[ T(u_t)^* T(u_t) ] \succeq \gamma I$.
We have the following minimax lower bound on the risk of any passive admissible algorithm $\calA$:
\begin{align}
  \inf_{\calA} \sup_{\theta \in \C^r} \E[\abs{\widehat{H}_{\calA} - \hinfnorm{H(g)}}] \geq C' \frac{\sigma}{M} \sqrt{\frac{r \log{r}}{N}} \:.
\end{align}
Here, $C'$ is a universal constant.
We note the form of the $\mathrm{poly}(r, M, 1/\gamma)$ can be recovered from the proof.
\end{mythm}

The significance of Theorem~\ref{thm:passive_lower_bound} is due to the fact that
under our query model, one can easily produce an estimate $\widehat{g} \in \C^r$ such that $\E[\hinfnorm{H(\widehat{g}) - H(g)}] \leq C'' \frac{\sigma}{M} \sqrt{\frac{r \log{r}}{N}}$,
for instance by setting $u_t = M e_1$ where $e_1$ is the first standard basis vector (see e.g. Theorem 1.1 of \cite{tuboczar}).
That is, for passive algorithms the number of samples to estimate the $r$ model parameters is equal (up to constants) to the number of samples needed to
estimate the $\hinf$-norm, at least for the worst-case.
The situation changes slightly when we look at active algorithms.
\begin{mythm}[Active lower bound]
\label{thm:active_lower_bound}
The following minimax lower bound on the risk of any active algorithm $\calA$ holds:
\begin{align}
  \inf_{\calA} \sup_{\theta \in \C^r} \E[\abs{\widehat{H}_{\calA} - \hinfnorm{H(g)}}] \geq C' \frac{\sigma}{M} \sqrt{\frac{r}{N}} \:.
\end{align}
Here, $C'$ is a universal constant.
\end{mythm}

We see in the active setting that the lower bound is weakened by a logarithmic
factor in $r$.  This bound shows that for low SNR regimes when $M  \ll \sqrt{r}$, the
gains of being active are minimal in the worst case.  On the other hand, for
high SNR regimes when $M  \gg \sqrt{r}$, the gains are potentially quite
substantial.  We are unaware of any algorithm that provably achieves the active
lower bound; it is currently unclear whether or not the lower bound is loose,
or if more careful design/analysis is needed to find a better active algorithm.
However, in Section~\ref{sec:active_lb_discussion} we discuss special cases of FIR filters
for which the lower bound in Theorem~\ref{thm:active_lower_bound} is not improvable.

We note that our proof draws heavily on techniques used to prove
minimax lower bounds for functional estimation, specifically in
estimating the $\ell_\infty$-norm of a unknown mean vector
in the Gaussian sequence model. An excellent overview of these techniques is given
in Chapter 22 of \cite{wu17}.

\subsection{Hypothesis testing and sector conditions}

An application that is closely related to estimating the $\hinf$-norm is
testing whether or not the $\hinf$-norm exceeds a certain fixed threshold
$\tau$. Specifically, consider a test statistic $\psi \in \{0, 1\}$ that discriminates between the two
alternatives $H_0 : \norm{H(g)}_{\hinf} \leq \tau$ and $H_1 : \norm{H(g)}_{\hinf} > \tau$.
This viewpoint is useful because it encompasses testing for more fine-grained characteristics of the
Nyquist plot $H(\omega)$ via simple transformations.
For instance, given $-\infty < a < 0 < b < +\infty$, one may test if
$H(\omega)$ is contained within a circle in the complex plane centered
at $(a+b)/2 + 0j$ with radius $(b-a)/2$ by equivalently checking if
the $\hinf$-norm of the system $H(g) - (a + b)/2$ is less than $(b-a)/2$;
this is known as the $[a,b]$-sector condition~\cite{gupta94sector}.

Due to the connection between estimation and hypothesis testing, our
results also give lower bounds on the sum of the Type-I and Type-II
errors of any test $\psi$ discriminating between the hypothesis $H_0$ and $H_1$.
Specifically, $\Omega( \frac{ \sigma^2 r \log{r}}{\tau^2 M} )$ queries
in the passive case (when $\tau$ is sufficiently small)
and $\Omega(\frac{\sigma^2 r}{\tau^2 M})$ queries in the active case are necessary
for any test to have Type-I and Type-II error less than a constant probability.

\subsection{Shortcomings of discretization}
\label{sec:active_lb_discussion}

In order to close the gap between the upper and lower bounds, one needs to
explicitly deal with the continuum of frequencies on $[0, 2\pi)$; here we argue that
if the maximizing frequency is known a-priori to lie in discrete set of points,
then the lower bound is sharp.

Suppose we consider a slightly different query model, where at each time $t$
the algorithm chooses a frequency $\omega_t \in [0, 2\pi)$ and receives
$y_t \sim \calN(H(\omega_t), 1)$.  For
simplicity let us also assume that the $\hinf$-norm is bounded by one.  Note
that by slightly enlarging $T(g)$ to the $2r \times 2r$ upper
left triangle, we can emulate this query model
by using a normalized complex sinusoid with
frequency $\omega_t$.

If the maximizing frequency for the $\hinf$-norm of the underlying system
is located on the grid $\{ 2\pi k / r \}_{k=0}^{r-1}$
and its phase is known,
this problem immediate reduces to a standard $r$-arm multi-armed bandit (MAB) problem
where each arm is associated with a point on the grid.
For this family of instances, the following active algorithm
has expected risk upper bounded by $\sqrt{r/N}$ times a constant:
\begin{enumerate}[(a)]
  \item Run a MAB algorithm that is optimal in the stochastic setting (such as MOSS \cite{audibert09}) for $N/2$ iterations,
    with each of the $r$ arms associated to a frequency $2\pi k / r$.
  \item Sample an index $I \in \{1, ..., r\}$ where:
    \begin{align*}
      \Pr(I = i) = \frac{ \text{\# number of times arm $i$ was pulled} }{N/2} \:.
    \end{align*}
  \item Query $\omega_I$ for $N/2$ times and return the sample mean.
\end{enumerate}
More generally, for any grid of frequencies $\{ \omega_k \}$ such that the
number of grid points is $O(r)$, this algorithm
obtains $O(\sqrt{r/N})$ expected risk. Hence the lower bound of Theorem~\ref{thm:active_lower_bound}
is actually sharp with respect to these instances.

Therefore, the issue that needs to be understood
is whether or not the continuum of frequencies on $[0, 2\pi)$
fundamentally requires additional sample complexity compared to a fixed discrete grid.
Note that a na{\"i}ve discretization argument is insufficient here.
For example, it is known (see e.g. Lemma 3.1 of \cite{tuboczar}) that
by choosing $P$ equispaced frequencies one obtains a discretization error
bounded by $O(r/P)$, e.g. $\norm{H(g)}_{\hinf} - \abs{H(\omega_k)} \leq O(r/P)$ for the largest $\omega_k$.
This bound is too weak, however,
since it requires that the number of arms scale as $O(r/\varepsilon)$ in order
to obtain a risk bounded by $\varepsilon$; in terms of $N$, the risk would scale $O(1/N^{1/3})$.

To summarize, if one wishes to improve the active lower bound
to match the rate given by Theorem~\ref{thm:passive_lower_bound}, one needs to consider
a prior distribution over hard instances where the support of the maximizing frequency is large (possibly infinite) compared to $r$.
On the other hand, if one wishes to construct an algorithm achieving the rate
of Theorem~\ref{thm:active_lower_bound}, then one will need to understand
the function $\omega \mapsto \abs{H(\omega)}$ at a much finer resolution than Lipschitz continuity.


\section{Proof of Main Results}

The proof of Theorem~\ref{thm:passive_lower_bound} and
Theorem~\ref{thm:active_lower_bound} both rely on a reduction to
Bayesian hypothesis testing. While this reduction is standard in the statistics
and machine learning communities (see e.g. Chapter 2 of \cite{tsybakov09}), we briefly
outline it here, as we believe these techniques are not as widely used in the controls literature.

First, let $\pi_1, \pi_2$ be two prior distributions on $\C^r$.
Suppose that for all $\theta_1 \in \pi_1$ we have $\hinfnorm{H(\theta_1)} = 0$
and for all $\theta_2 \in \pi_2$ we have $\hinfnorm{H(\theta_2)} = 2c$
for some $c > 0$.
Let $\Pr_{\pi_i}$ denote the joint distribution of $( u_1, y_1, ..., u_N, y_N)$,
which combines the prior distribution $\pi_i$ with the observation model.
Then,
\begin{align*}
  \sup_{\theta \in \C^r} \E[\abs{ \widehat{H} - \hinfnorm{H(\theta)}}] &\geq c \sup_{\theta \in \C^r} \Pr_\theta( \abs{\widehat{H} - \hinfnorm{H(\theta)}} \geq c ) \\
  &\geq c \max_{i=1,2}\left\{ \int \Pr_\theta( \abs{\widehat{H} - \hinfnorm{H(\theta)}} \geq c ) \; \pi_i(d\theta) \right\} \\
  &\geq \frac{c}{2}\left( \Pr_{\pi_1}( \abs{\widehat{H}} \geq c ) + \Pr_{\pi_2}(\abs{\widehat{H} - 2c} \geq c ) \right) \\
  &\geq \frac{c}{2}\left( \Pr_{\pi_1}( \abs{\widehat{H}} \geq c ) + \Pr_{\pi_2}(\abs{\widehat{H}} < c ) \right) \\
  &\geq \frac{c}{2}(1 - \dTV(\Pr_{\pi_1}, \Pr_{\pi_2})) \:,
\end{align*}
where for two measures $\Pr, \bbQ$ we define the total-variation (TV) distance
as $\dTV(\Pr, \bbQ) = \sup_{A} \abs{\Pr(A) - \bbQ(A)}$.
Hence, if one can construct two prior distributions $\pi_1, \pi_2$ with the
aforementioned properties and furthermore show that $\dTV(\Pr_{\pi_1}, \Pr_{\pi_2}) \leq 1/2$,
then one deduces that the minimax risk is lower bounded by $c/4$.
This technique is generally known as Le Cam's method, and will be our high-level proof strategy.

As working directly with the TV distance is often intractable, one typically computes
upper bounds to the TV distance. We choose to work with both the KL-divergence and the $\chi^2$-divergence.
The KL-divergence is defined as $\dKL(\Pr, \bbQ) = \int \log\left(\frac{d\Pr}{d\bbQ}\right) \: d\Pr$,
and the $\chi^2$-divergence is defined as $\dChiSq(\Pr, \bbQ) = \int \left( \frac{d\Pr}{d\bbQ} - 1\right)^2 \; d\bbQ$
(we assume that $\Pr \ll \bbQ$ so these quantities are well-defined).
One has the standard inequalities
$\dTV(\Pr, \bbQ) \leq \sqrt{\frac{1}{2} \dKL(\Pr, \bbQ)}$ and $\dTV(\Pr, \bbQ) \leq \sqrt{\dChiSq(\Pr, \bbQ)}$
\cite{tsybakov09}.

\subsection{Proof of passive lower bound (Theorem~\ref{thm:passive_lower_bound})}

The main reason for working with the $\chi^2$-divergence is that it operates nicely with mixture distributions,
as illustrated by the following lemma.
\begin{mylemma}[see e.g. Lemma 22.1 of \cite{wu17}]
Let $\Theta$ be a parameter space and for each $\theta \in \Theta$ let $\Pr_\theta$ be a measure over $\calX$ indexed by $\theta$. Fix a measure $\bbQ$ on $\calX$ and a prior measure $\pi$ on $\Theta$.
Define the mixture measure $\Pr_\pi = \int \Pr_\theta \; \pi(d\theta)$.
Suppose for every $\theta \in \Theta$, the measures $\Pr_\theta$ and $\bbQ$ are both absolutely continuous w.r.t. a fixed base measure $\mu$ on $\calX$.
Define the function $G(\theta_1, \theta_2)$ as
\begin{align*}
  G(\theta_1, \theta_2) := \int \frac{ \frac{d\Pr_{\theta_1}}{d\mu}  \frac{d\Pr_{\theta_2}}{d\mu}   }{  \frac{d\bbQ}{d\mu}   } \; \mu(dx) \:.
\end{align*}
We have that:
\begin{align*}
  \dChiSq(\Pr_\pi , \bbQ) = \E_{\theta_1, \theta_2 \sim \pi^{\otimes 2}}[ G(\theta_1, \theta_2) ] - 1 \:.
\end{align*}
\end{mylemma}

We now specialize this lemma to our setting. Here, our distributions $\Pr_\theta$ are over
$(u_1, x_1, ..., u_N, x_N)$;
for a fixed system parameter $\theta \in \C^r$, the joint distribution $\Pr_\theta$
has the density (assuming that $u_t$ has the density $\gamma_t(u_t)$):
\begin{align*}
    p_\theta(u_1, x_1, ..., u_N, x_N) = \prod_{t=1}^{N} \gamma_t(u_t) \phi(x_t; T(\theta) u_t) \:,
\end{align*}
where $\phi(x; \mu)$ denotes the PDF of the multivariate Gaussian $\calN( \mu, \sigma^2 I)$.
Note that this factorization with $\gamma_t(\cdot)$ independent of $\theta$ is only possible
under the passive assumption.

\begin{mylemma}
\label{lemma:g_one_two}
Supposing that $u_t \sim \gamma_t(\cdot)$, we have that
\begin{align*}
  G(\theta_1, \theta_2) = \E_{u_t}\left[ \exp\left(\frac{1}{\sigma^2} \sum_{t=1}^{N}\Re(\ip{T(\theta_1)u_t}{T(\theta_2)u_t}) \right) \right].
\end{align*}
\end{mylemma}
\begin{proof}
We write:
\begin{align*}
  G(\theta_1, \theta_2) &= \int \frac{p_{\theta_1} p_{\theta_2}}{p_0} \; du_1 dx_1 ... du_N dx_N \\
    &= \int \prod_{t=1}^{N} \frac{\gamma_t(u_t) \phi(x_t; T(\theta_1) u_t) \phi(x_t; T(\theta_2) u_t)}{\phi(x_t; 0)} \; du_t dx_t \\
    &= \prod_{t=1}^{N} \int \frac{\gamma_t(u_t) \phi(x_t; T(\theta_1) u_t) \phi(x_t; T(\theta_2) u_t)}{\phi(x_t; 0)} \; du_t dx_t \\
    &\stackrel{(a)}= \prod_{t=1}^{N} \E_{u_t}\left[ \exp\left(\frac{1}{\sigma^2} \Re(\ip{T(\theta_1)u_t}{T(\theta_2)u_t}) \right) \right] \\
    &\stackrel{(b)}{=} \E_{u_t}\left[ \exp\left(\frac{1}{\sigma^2} \sum_{t=1}^{N}\Re(\ip{T(\theta_1)u_t}{T(\theta_2)u_t}) \right) \right] \:.
\end{align*}
In part (a) we complete the square
and in part (b) we use the fact that the distributions for $u_t$ are independent as a consequence of the passive assumption.
In particular for (a), we first observe that:
\begin{align*}
  &\norm{x_t - T(\theta_1)u_t}_2^2 + \norm{x_t - T(\theta_2)u_t}_2^2 - \norm{x_t}_2^2 \\
  &\qquad= 2\norm{x_t}_2^2 + \norm{T(\theta_1)u_t}_2^2 + \norm{T(\theta_2)u_t}_2^2 - 2\Re(\ip{x_t}{T(\theta_1) u_t + T(\theta_2)}u_t) - \norm{x_t}_2^2 \\
  &\qquad= \norm{x_t}_2^2 + \norm{T(\theta_1)u_t + T(\theta_2)u_t}_2^2 - 2\Re(\ip{x_t}{T(\theta_1) u_t + T(\theta_2)}u_t) - 2\Re(\ip{T(\theta_1)u_t}{T(\theta_2)u_t}) \\
  &\qquad= - 2\Re(\ip{T(\theta_1)u_t}{T(\theta_2)u_t}) + \norm{x_t - (T(\theta_1)u_t + T(\theta_2)u_t)}_2^2 \:.
\end{align*}
Hence writing $\phi(x; \mu) = C \exp(-\frac{1}{2\sigma^2} \norm{x-\mu}_2^2)$, we obtain:
\begin{align*}
  &\int \frac{\phi(x_t; T(\theta_1) u_t) \phi(x_t; T(\theta_2) u_t)}{\phi(x_t; 0)} \; dx_t = \int C \exp\left( -\frac{1}{2\sigma^2} ( \norm{x_t - T(\theta_1)u_t}_2^2 + \norm{x_t - T(\theta_2)u_t}_2^2 - \norm{x_t}_2^2  )  \right) \; dx_t \\
  &\qquad= \exp\left(\frac{1}{\sigma^2} \Re(\ip{T(\theta_1)u_t}{T(\theta_2)u_t})\right) \int C \exp\left( -\frac{1}{2\sigma^2} \left( \norm{x_t - (T(\theta_1)u_t + T(\theta_2)u_t)}_2^2 \right) \right) dx_t \\
  &\qquad= \exp\left(\frac{1}{\sigma^2} \Re(\ip{T(\theta_1)u_t}{T(\theta_2)u_t})\right) \int \phi(x_t; T(\theta_1) u_t + T(\theta_2) u_t) \; dx_t \\
  &\qquad= \exp\left(\frac{1}{\sigma^2} \Re(\ip{T(\theta_1)u_t}{T(\theta_2)u_t})\right) \:.
\end{align*}
\end{proof}

We now construct two prior distributions on $\theta$.
The first prior will be the system with all coefficients zeros, i.e. $\pi_1 = \{ 0 \}$. The second prior
will be more involved.
To construct it, we let $\Sigma := \frac{1}{N} \sum_{t=1}^{N} \E_{u_t \sim \gamma_t}[ T(u_t)^* T(u_t) ]$.
By the admissibility assumption on the algorithm $\calA$, $\Sigma$ is invertible.
Let $\calI \subseteq \{1, ..., r\}$ denote an index set to be specified.
Let $F \in \C^{r \times r}$ denote the unnormalized discrete Fourier transform (DFT) matrix
(i.e. $FF^* = rI$ and $F^{-1} = \frac{1}{r} F^*$).
We define our prior distribution $\pi_2$ as, for some $\tau > 0$ to be chosen:
\begin{align*}
  \pi_2 = \mathrm{Unif}(\{ \tau \Sigma^{-1/2} F^{-1} e_i \}_{i \in \calI}) \:.
\end{align*}
We choose $\calI$ as according to the following proposition.
\begin{myprop}
\label{prop:index_set}
Let $u_1, ..., u_N \in \C^r$ be independently drawn from $N$ distributions such that $\norm{u_t}_2 \leq M$ a.s for all $t = 1 , ..., N$.
Let $\Sigma = \frac{1}{N} \sum_{t=1}^{N} \E_{u_t}[ T(u_t)^* T(u_t) ]$
and suppose that $\Sigma$ is invertible.
There exists an index set $\calI \subseteq \{1, ..., r\}$ such that $\abs{\calI} \geq r/2$ and
for all $i \in \calI$,
\begin{align*}
  \norm{ F \Sigma^{-1/2} F^{-1} e_i }_\infty \geq \frac{1}{2M} \:.
\end{align*}
\end{myprop}
\begin{proof}
First, we observe that:
\begin{align*}
  \sum_{i=1}^{r} e_i^\T F \Sigma^{1/2} F^{-1} e_i = \Tr(F \Sigma^{1/2} F^{-1}) = \Tr(\Sigma^{1/2}) \:.
\end{align*}
Let $\lambda_1, ..., \lambda_r$ denote the eigenvalues of $\Sigma$. By Cauchy-Schwarz,
\begin{align*}
  \Tr(\Sigma^{1/2}) = \sum_{i=1}^{r} \sqrt{\lambda_i} \leq \sqrt{r \Tr(\Sigma)} \:.
\end{align*}
Now for any fixed $u \in \C^r$ satisfying $\norm{u}_2 \leq M$, we have:
\begin{align*}
  \Tr( T(u)^* T(u) ) \leq r M^2 \:.
\end{align*}
Hence,
\begin{align*}
  \Tr(\Sigma) = \frac{1}{N} \sum_{t=1}^{N} \E_{u_t}[ \Tr(T(u_t)^* T(u_t)) ] \leq r M^2 \:.
\end{align*}
That is, we have shown that:
\begin{align*}
  \sum_{i=1}^{r} e_i^\T F \Sigma^{1/2} F^{-1} e_i \leq r M \Longleftrightarrow \frac{1}{r} \sum_{i=1}^{r} e_i^\T F \Sigma^{1/2} F^{-1} e_i \leq M \:.
\end{align*}
Now we state an auxiliary proposition whose proof follows from Markov's inequality.
\begin{myprop}
Let $a_1, ..., a_r \in \R$ satisfy $a_i \geq 0$ and $\frac{1}{r} \sum_{i=1}^{r} a_i \leq M$.
Then there exists an index set $\calI \subseteq \{1, ..., r\}$ with cardinality $\abs{\calI} \geq r/2$ such that
$a_i \leq 2M$ for all $i \in \calI$.
\end{myprop}
By the auxiliary proposition, there exists an index set $\calI$ such that
$\abs{\calI} \geq r/2$ and $e_i^\T F \Sigma^{1/2} F^{-1} e_i \leq 2 M$ for all $i \in \calI$.
Hence, for any $i \in \calI$,
\begin{align*}
  \norm{ F \Sigma^{-1/2} F^{-1} e_i }_\infty &\geq e_i^\T F \Sigma^{-1/2} F^{-1} e_i \geq \frac{1}{e_i^\T F \Sigma^{1/2} F^{-1} e_i} \geq \frac{1}{2M} \:.
\end{align*}
Above, the first inequality holds because the $\ell_\infty$-norm of the $i$-th column of a matrix exceeds the absolute value of the $i, i$-th position of the matrix, the second
inequality is because for any positive definite matrix $M$, we have $(M^{-1})_{ii} \geq 1/M_{ii}$
and the last inequality is due to the property of $\calI$.
\end{proof}

We now observe that for indices $\ell_1, \ell_2 \in \calI$, defining $\Delta := \sum_{t=1}^{N} (T(u_t)^* T(u_t) - \E_{u_t}[T(u_t)^* T(u_t)])$:
\begin{align*}
    \sum_{t=1}^{N} \ip{T( \Sigma^{-1/2} F^{-1} e_{\ell_1}) u_t }{ T( \Sigma^{-1/2} F^{-1} e_{\ell_2}) u_t  } &\stackrel{(a)}{=} \sum_{t=1}^{N} \ip{T(u_t) \Sigma^{-1/2} F^{-1} e_{\ell_1}}{ T(u_t) \Sigma^{-1/2} F^{-1} e_{\ell_2}} \\
    &= \sum_{t=1}^{N} e_{\ell_1}^\T F^{-*} \Sigma^{-1/2} T(u_t)^* T(u_t) \Sigma^{-1/2} F^{-1} e_{\ell_2} \\
    &= e_{\ell_1}^\T F^{-*} \Sigma^{-1/2} ( N\Sigma + \Delta) \Sigma^{-1/2} F^{-1} e_{\ell_2} \\
    &= N e_{\ell_1}^\T F^{-*} F^{-1} e_{\ell_2} + e_{\ell_1}^* F^{-*} \widetilde \Delta F^{-1} e_{\ell_2} \\
    &= \frac{N}{r} \ind_{\ell_1=\ell_2} + e_{\ell_1}^\T F^{-*} \widetilde \Delta F^{-1} e_{\ell_2} \:,
\end{align*}
where $\widetilde \Delta:=\Sigma^{-1/2} \Delta \Sigma^{-1/2}$. In (a) we used the fact that for two vectors $u, v$ we have $T(u) v = T(v) u$,
e.g. the fact that convolution is commutatitive.
Combining this calculation with Lemma~\ref{lemma:g_one_two},
\begin{align*}
  \E_{\theta_i}[ G(\theta_1, \theta_2) ] &= \E_{\ell_i, u_t}\left[ \exp\left( \frac{\tau^2}{\sigma^2}  \frac{N}{r} \ind_{\ell_1=\ell_2} \right) \exp\left( \frac{\tau^2}{\sigma^2} \Re\left(e_{\ell_1}^\T F^{-*} \widetilde \Delta F^{-1} e_{\ell_2}  \right) \right) \right] \\
    &\stackrel{(a)}{\leq} \sqrt{\E_{\ell_i, u_t}\left[ \exp\left(\frac{2 \tau^2 N}{\sigma^2 r} \ind_{\ell_1 = \ell_2}\right) \right]} \sqrt{\E_{\ell_i, u_t}\left[ \exp\left(\frac{2\tau^2}{\sigma^2} \Re(e_{\ell_1}^\T F^{-*} \widetilde \Delta F^{-1} e_{\ell_2}) \right)\right]} \\
    &\stackrel{(b)}{\leq} \sqrt{ \exp\left(\frac{2 N \tau^2}{\sigma^2 r}\right) \frac{2}{r} + 1 - \frac{2}{r}} \sqrt{\E_{\ell_i, u_t}\left[ \exp\left(\frac{2\tau^2}{\sigma^2} \Re(e_{\ell_1}^\T F^{-*} \widetilde \Delta F^{-1} e_{\ell_2}) \right)\right]} \:.
\end{align*}
where in (a) we used Cauchy-Schwarz
and in (b) we used the fact that $\abs{\calI} \geq r/2$.
Now condition on $\ell_1, \ell_2$.
For a $1 \leq t \leq N$, define the random variable $\psi_t$ as:
\begin{align*}
  \psi_t &:= \Re( e_{\ell_1}^\T F^{-*} \Sigma^{-1/2} T(u_t)^* T(u_t) \Sigma^{-1/2} F^{-1} e_{\ell_2} ) \\
  &\qquad - \Re( e_{\ell_1}^\T F^{-*} \Sigma^{-1/2} \E_{u_t}[T(u_t)^* T(u_t)] \Sigma^{-1/2} F^{-1} e_{\ell_2} ) \:.
\end{align*}
We have that $\E_{u_t}[\psi_t] = 0$ by construction.
Furthermore, note that $\norm{F^{-1} e_{\ell}}_2 = 1/\sqrt{r}$
for $\ell=1, ..., r$ and also that that $\norm{T(u)} \leq \norm{H(u)}_{\hinf} \leq \norm{u}_1 \leq \sqrt{r} \norm{u}_2$
for any vector $u \in \C^r$. These facts, along with the assumption that $\Sigma \succeq \gamma I$, show that $\abs{\psi_t} \leq 2 M^2 / \gamma$ almost surely.
Hence, $\sum_{t=1}^{N} \psi_t$ is a zero-mean sub-Gaussian random variable with sub-Gaussian parameter
$4 M^4 N / \gamma^2$ (see e.g. Ch. 2 of \cite{vershynin18} for background exposition on sub-Gaussian random variables).
Therefore, we know that for any $t > 0$, its moment generating function (MGF) is bounded as:
\begin{align*}
  \E_{u_t | \ell_i}\left[ \exp\left( t \sum_{t=1}^{N} \psi_t \right)\right] \leq \exp( 2 t^2 M^4 N / \gamma^2 ) \:.
\end{align*}
Hence by iterating expectations and setting $t = 2\tau^2/\sigma^2$, we have:
\begin{align*}
  \E_{\ell_i, u_t}\left[ \exp\left(\frac{2\tau^2}{\sigma^2} \Re(e_{\ell_1}^\T F^{-*} \Sigma^{-1/2} \Delta \Sigma^{-1/2} F^{-1} e_{\ell_2}) \right)\right] \leq \exp( 8 \tau^4 M^4 N / (\sigma^4 \gamma^2  )) \:.
\end{align*}
Therefore, for any choice of $\tau$ such that:
\begin{align}
  \tau^4 \leq (\log(1.1)/8) \sigma^4 \gamma^2 / (M^4 N) \:, \label{eq:condone}
\end{align}
we have that:
\begin{align*}
  \E_{\theta_i, u_t}[ G(\theta_1, \theta_2) ] \leq \sqrt{1.1}\sqrt{ \exp\left(\frac{2 N \tau^2}{\sigma^2 r}\right) \frac{2}{r} + 1 - \frac{2}{r}} \:.
\end{align*}
Hence if $r \geq 5$ and if we set $\tau$ to be:
\begin{align}
  \tau^2 = \frac{\sigma^2 r \log(0.211 r)}{2N} \:, \label{eq:tausetting}
\end{align}
we have that:
\begin{align*}
  \dChiSq(\Pr_\pi, \Pr_0) = \E_{\theta_i, u_t}[ G(\theta_1, \theta_2) ] - 1 \leq 1/4 \:,
\end{align*}
assuming the condition \eqref{eq:condone} is satisfied.
This bound then implies that $\dTV(\Pr_\pi, \Pr_0) \leq 1/2$.
We now aim to choose $N$ so that the condition \eqref{eq:condone} holds.
Plugging our setting of $\tau$ from \eqref{eq:tausetting} into \eqref{eq:condone} and rearranging yields the condition
$N \geq (2/\log(1.1)) r^2 \log^2(0.211 r) M^4 / \gamma^2 = \widetilde{\Omega}( r^2 M^4 / \gamma^2)$.

To conclude, we need to show a minimum separation between
the $\hinf$-norm on $\pi_1$ vs. $\pi_2$. Clearly $\hinfnorm{H(\theta)} = 0$ on $\pi_1$.
On the other hand, for $\theta \in \pi_2$, we observe that for $i \in \calI$,
\begin{align*}
  \norm{H(\tau \Sigma^{-1/2} F^{-1} e_i)}_{\hinf} &\stackrel{(a)}{\geq} \tau \norm{F \Sigma^{-1/2} F^{-1} e_i}_\infty \stackrel{(b)}{\geq} \frac{\tau}{2 M} \:,
\end{align*}
where inequality (a) comes from $\hinfnorm{H(g)} \geq \norm{F g}_\infty$ for any $g$ and
inequality (b) comes from Proposition~\ref{prop:index_set}.
Hence we have constructed two prior distributions with a separation of $c = \Omega(\tau/M)$
but a total variation distance less than $1/2$. Theorem~\ref{thm:passive_lower_bound} now follows.

\subsection{Proof of active lower bound (Theorem~\ref{thm:active_lower_bound})}

For this setting we let $\pi_1 = \{ 0 \}$ and
$\pi_2 = \{ \tau F^{-1} e_i \}_{i=1}^{r}$.
The proof proceeds by bounding the KL-divergence $\dKL(\Pr_{\pi_1}, \Pr_{\pi_2})$.
To do this, we first bound $\dKL(\Pr_{0}, \Pr_{i})$, where
$\Pr_{0}$ is the joint distribution induced by the parameter $g = 0$ and
$\Pr_{i}$ is the joint distribution induced by the parameter $g = \tau F^{-1} e_i$.
Proceeding similarly to the proof of Theorem 1.3 in \cite{tuboczar},
\begin{align*}
  \dKL(\Pr_{0}, \Pr_{i}) &= \E_{\Pr_{0}}\left[ \log \prod_{t=1}^{N} \frac{\gamma_t(u_t | \{u_{k}, x_{k}\}_{k=1}^{t-1} ) p_0(x_t | u_t) }{\gamma_t(u_t | \{u_{k}, x_{k}\}_{k=1}^{t-1} ) p_i(x_t | u_t) } \right] \\
    &= \E_{\Pr_{0}}\left[ \log \prod_{t=1}^{N} \frac{p_0(x_t | u_t) }{p_i(x_t | u_t) } \right] \\
    &= \sum_{t=1}^{N} \E_{\Pr_{0}}\left[ \log\frac{p_0(x_t | u_t) }{p_i(x_t | u_t) }  \right] \\
    &= \sum_{t=1}^{N} \E_{u_t \sim \Pr_{0}}[ \dKL(\calN(0, \sigma^2 I), \calN( T(\tau F^{-1} e_i)u_t, \sigma^2 I)) ] \\
    &= \frac{\tau^2}{2\sigma^2} \sum_{t=1}^{N} \E_{u_t \sim \Pr_{0}}[ \norm{T( F^{-1} e_i) u_t}_2^2 ] \:.
\end{align*}
A simple calculation shows that:
\begin{align*}
    \sum_{i=1}^{r} T(F^{-1} e_i)^* T(F^{-1} e_i) = \diag\left(1, \frac{r-1}{r}, \frac{r-2}{r}, ..., \frac{1}{r}\right) \:,
\end{align*}
and hence the operator norm of this matrix is bounded by one.
Therefore, by convexity of $\dKL$,
\begin{align*}
  \dKL(\Pr_{\pi_1}, \Pr_{\pi_2}) &\leq \frac{1}{r} \sum_{i=1}^{r} \dKL(\Pr_0, \Pr_i) \\
    &= \frac{\tau^2}{2 \sigma^2 r} \sum_{i=1}^{r} \sum_{t=1}^{N} \E_{u_t \sim \Pr_{0}}[ \norm{T( F^{-1} e_i) u_t}_2^2 ] \\
    &= \frac{\tau^2}{2 \sigma^2 r} \sum_{t=1}^{N} \sum_{i=1}^{r} \E_{u_t \sim \Pr_{0}}[ \norm{T( F^{-1} e_i) u_t}_2^2 ] \\
    &= \frac{\tau^2}{2 \sigma^2 r} \sum_{t=1}^{N} \E_{u_t \sim \Pr_{0}}\left[u_t^* \left( \sum_{i=1}^{r} T(F^{-1} e_i)^* T(F^{-1} e_i) \right) u_t \right] \\
    &\leq \frac{\tau^2 N M^2}{2 \sigma^2 r}  \bignorm{\sum_{i=1}^{r} T(F^{-1} e_i)^* T(F^{-1} e_i)} \\
    &\leq \frac{\tau^2 N M^2}{2 \sigma^2 r} \:.
\end{align*}
Hence if we set $\tau = \frac{\sigma}{M} \sqrt{\frac{r}{N}}$, we have by Pinsker's inequality that
$\dTV(\Pr_{\pi_1}, \Pr_{\pi_2}) \leq 1/2$. Finally, we note that $\hinfnorm{H(\tau F^{-1}e_i)} \geq \tau$ and conclude.


\section{Experiments}

We conduct experiments comparing a simple non-adaptive
estimator based on least-squares (which we call the \emph{plugin} estimator)
to three active algorithms: two similar
algorithms essentially based on the power method
\cite{rojas2012analyzing,wahlberg2010non} and one based on weighted Thompson
Sampling (WTS) \cite{muller17bandit}. Pseudocode for the plugin estimator is shown
in Algorithm~\ref{alg:plugin}.
For completeness, in the appendix we describe the power method based algorithms in Algorithms \ref{alg:power_method_A} and \ref{alg:power_method_B}, and
the WTS algorithm in Algorithm~\ref{alg:WTS}.
For brevity, we assume input normalization of $\|u^{(t)}\|_2 = 1$; for different SNR the algorithms are modified accordingly.
\begin{algorithm}[htb]
\caption{Plugin Estimator}
\label{alg:plugin}
Input: Normalized $\{u^{(t)}\}$.\\
\For{$t = 1$ \upshape to $N$}
{
   Perform the experiment $y^{(t)} = Gu^{(t)} + \eta^{(t)}$.
}
Form $\hat G$ from a least-squares fit of $\{y^{(t)}\}$ and $\{u^{(t)}\}$.\\
\Return $\hat H = \|\hat G\|_{\hinf}$.
\end{algorithm}

We compare the performance of these four algorithms on a suite of
random plants and random draws of noise.
We note, however, that it is difficult to place these algorithms on even footing when making a comparison, especially in the presence of output noise. Reasons for this are:
 \begin{itemize}
   \item The parameters deemed ``fixed'' may be beneficial (or adversarial) to one algorithm or another.
   \item The amount of ``side information'' (e.g. noise covariance) an algorithm expects to receive may not be comparable across algorithms.
 \end{itemize}
For an example of the first point, a large experiment budget is beneficial to
the plugin and WTS estimators as they generally obtain better estimates with
each new experiment while power method estimators hit a ``noise floor'' and
stop improving.

The plants we test are of the form $G(z) = \; \sum_{k=0}^{r-1}z^{-k}\rho^{k}\eta_k
$, where $\rho\in(0,1]$ and $\eta_k \stackrel{\text{i.i.d}}{\sim}\text{Unif}[-1,1]$. For each suite of tests, we hold all other parameters fixed as shown in Table \ref{par_table}.
To compare the aggregate performance across suites of random plants, we use \emph{performance profiles}, a tool in the optimization community popularized by Dolan and Mor\'e \cite{dolan2002benchmarking}. Given a suite of methods $\{m_i\}$ and a metric $d(m_i,m_j)$ for per-instance performance, performance profiles show the percentage of instances where a particular method $m$ is within $\tau$ of the best method. In our case, the metric will be the difference in relative error between the method's estimate and the true $\Hinf$-norm of the plant under consideration. As an example, \textsc{Plugin}(.05) would be the percentage of instances where the relative error of the plugin estimator is within $5$ percentage points of the smallest relative error on that instance. Performance profiles are meant to show broad differences in algorithm performance and are robust to per-instance variation in the data, when interpreted correctly.
\begin{figure*}[htb]
\centering
\subfloat[High SNR, decay]{\includegraphics[width=3in]{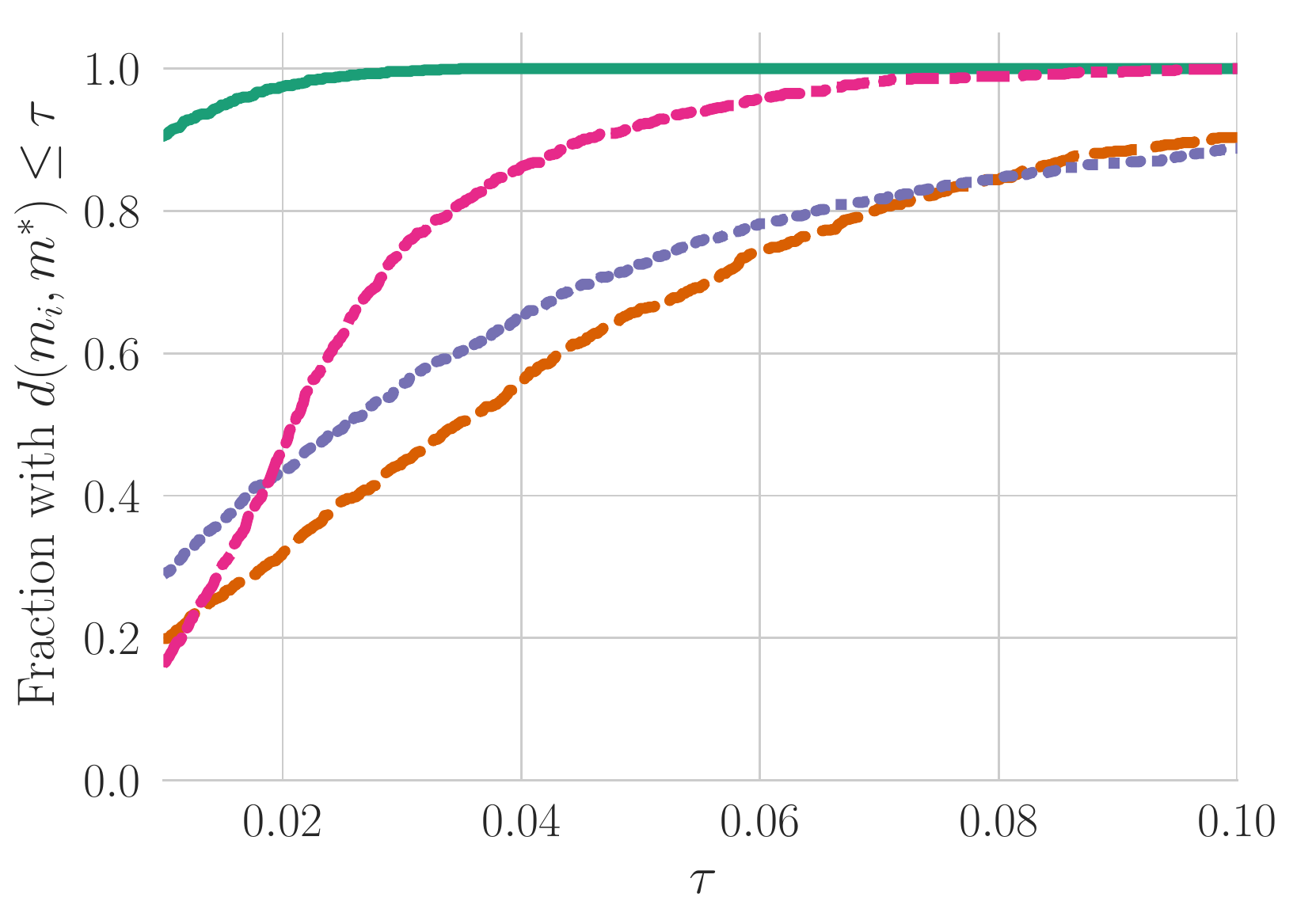}
\label{fig:hi_dec}}
\hfil
\subfloat[Low SNR, decay]{\includegraphics[width=3in]{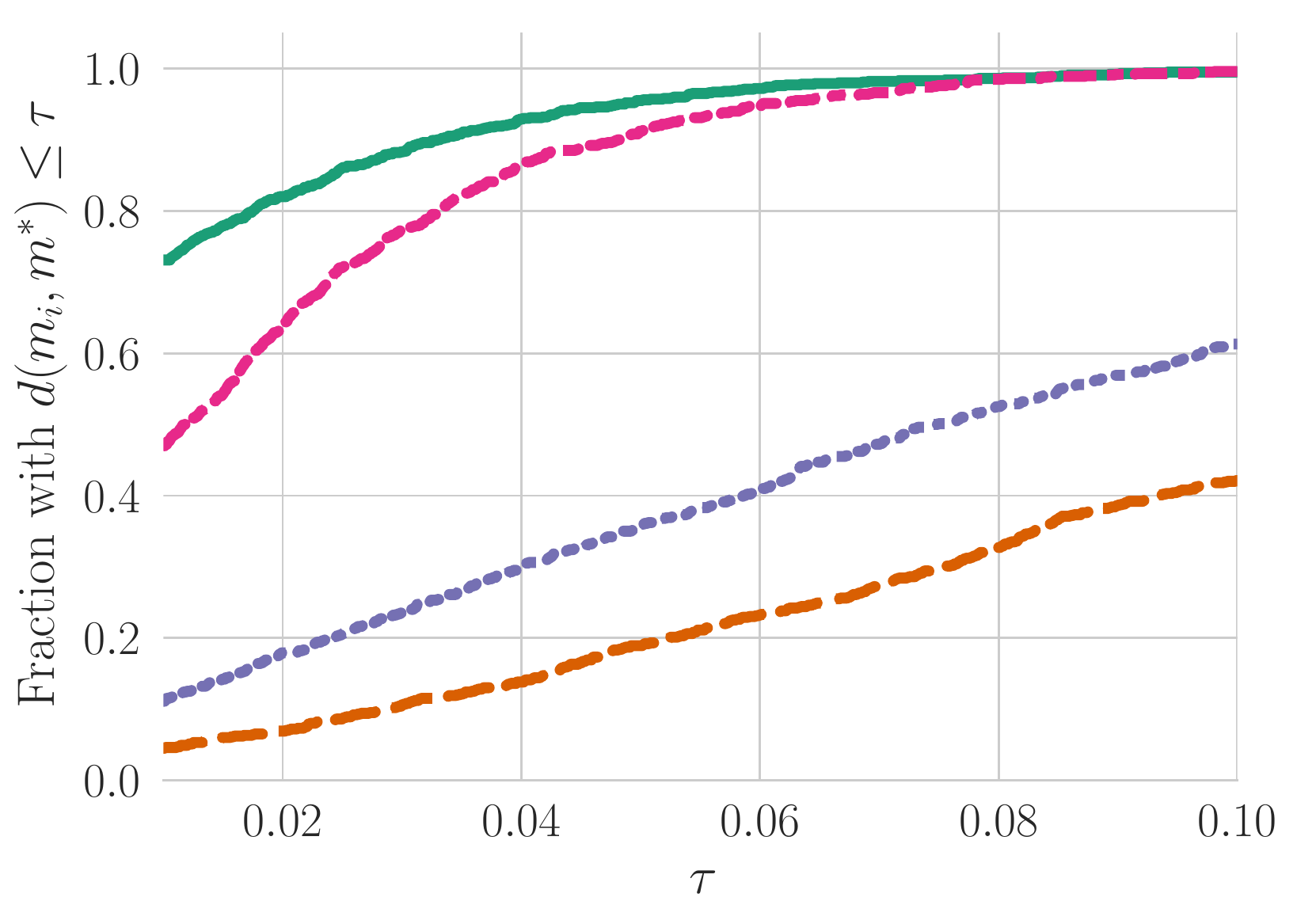}
\label{fig:lo_dec}}\\
\subfloat[High SNR, no decay]{\includegraphics[width=3in]{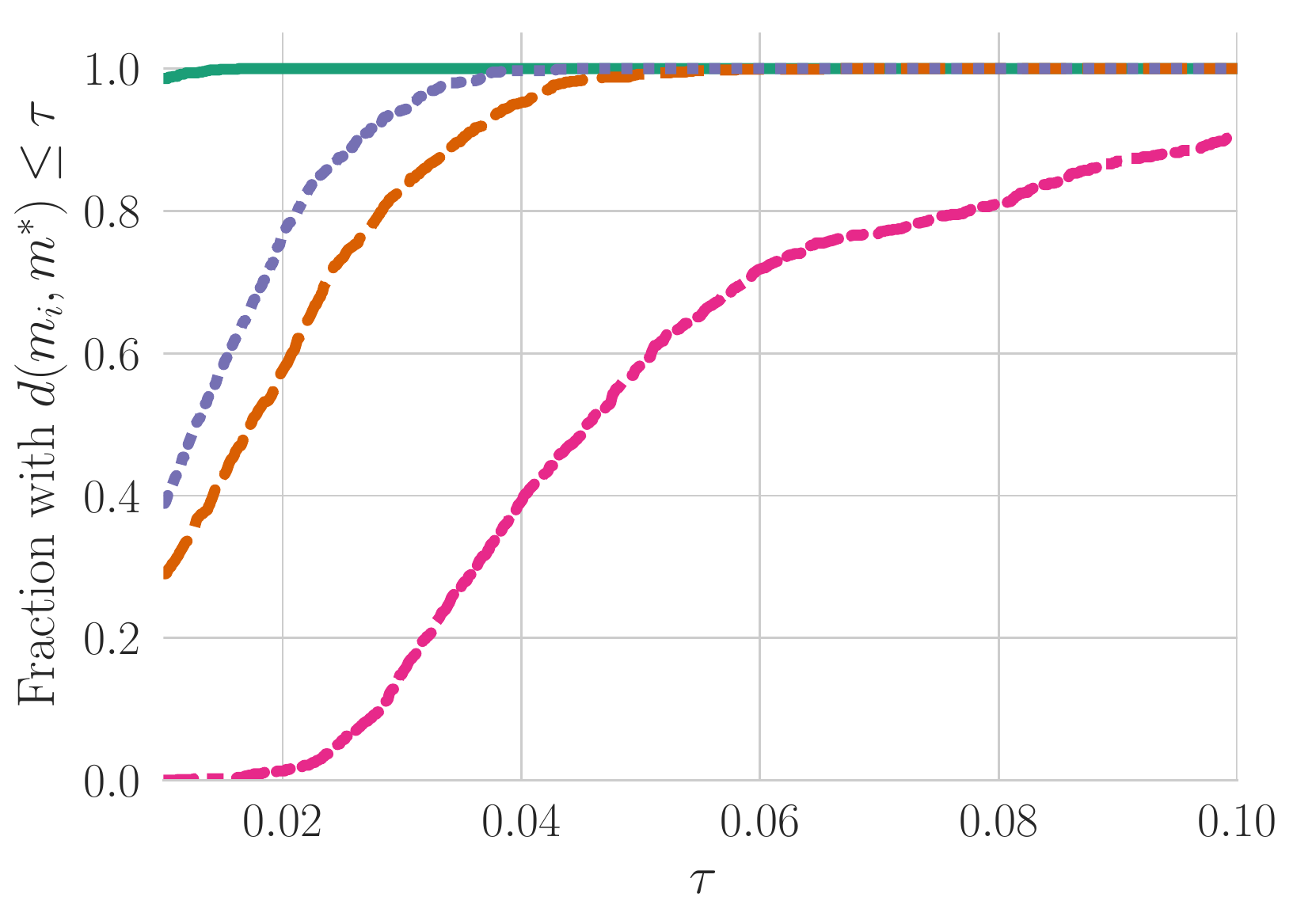}
\label{fig:hi_nodec}}
\hfil
\subfloat[Low SNR, no decay]{\includegraphics[width=3in]{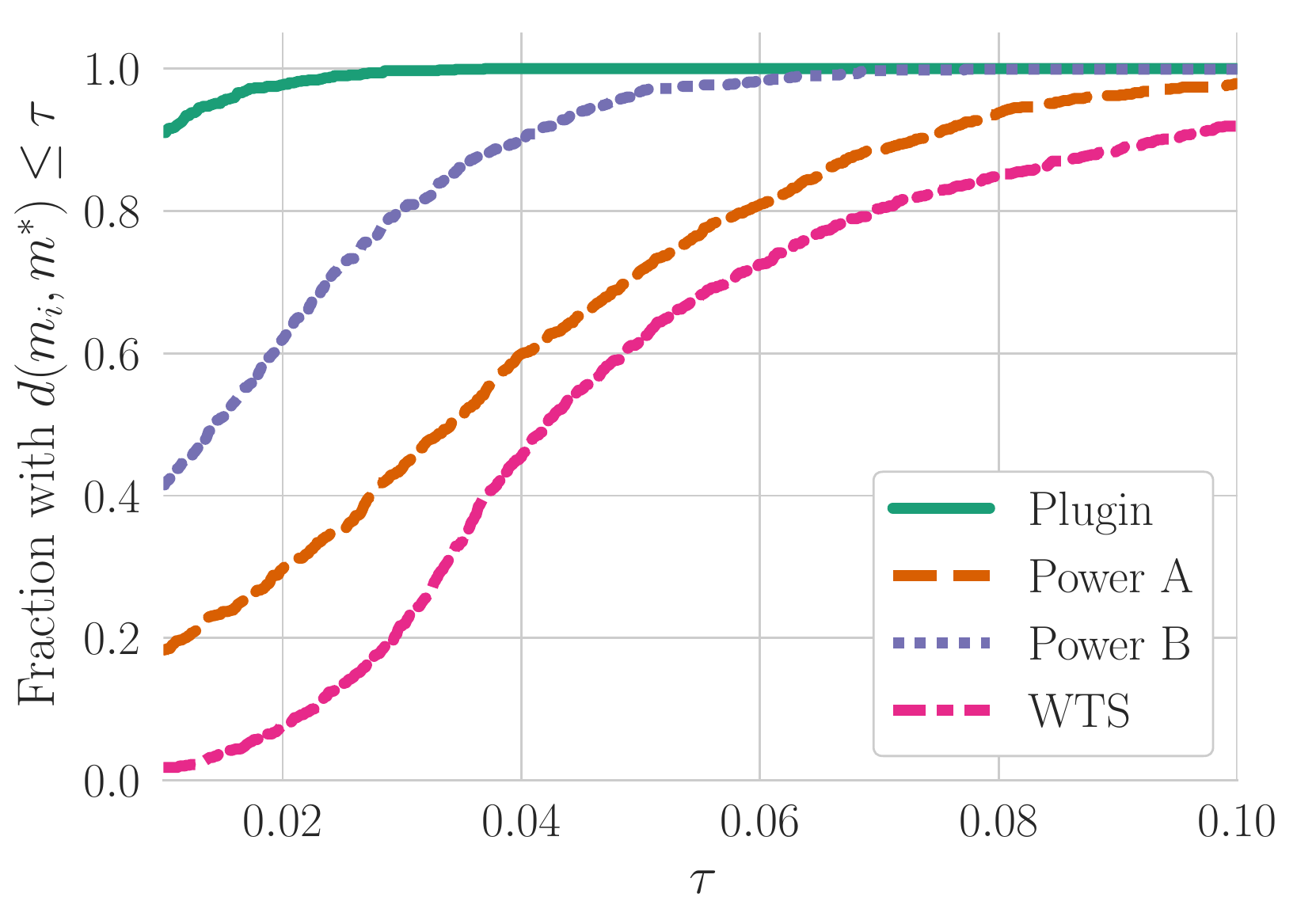}
\label{fig:lo_nodec}}
\caption{Performance profiles for the plugin, power method A, power method B, and weighted Thompson Sampling estimators.}
\label{fig:all}
\end{figure*}
\begin{table}[!h]
\footnotesize
\renewcommand{\arraystretch}{1.5}
\caption{Experiment Parameters}
\centering
\begin{tabular}{cc}
\hline
\bfseries Type & Value\\
\hline
SNR $\frac{\|u\|_2}{\sigma}$ & 20 (high), 10 (low)\\
Experiment budget $N$ & 200\\
Plant length $r$ & 10 \\
Input/output data length $r'$ & 50 \\
\hline
\hline
Plant decay $\rho$ &  0.75 (decay), 1.0 (no decay)\\
Number of random plants & 100 \\
Noise instances per plant & 10
\end{tabular}\label{par_table}
\end{table}
With this in mind, the performance profiles comparing the four algorithms\footnote{The code producing these plots can be found at \texttt{https://github.com/rjboczar/hinf-lower-bounds-ACC}. The experiments were carried out using the PyWren framework\cite{jonas2017occupy}.} are shown in Fig.~\ref{fig:all}. We see that for the plants with decaying impulse response coefficients, the plugin and WTS estimators are comparable. The latter estimator performs relatively worse in the experiments corresponding to no coefficient decay. As alluded to previously, this can most likely be attributed to our particular experimental setup: the WTS algorithm effectively grids the frequency response curve of the plant, and this set of plants allows for more relative variation in the curve than the onces with decaying coefficients.

\section{Conclusion}

We study lower bounds for $\hinf$-norm estimation for both passive and active
algorithms. Our analysis shows that in the passive case model identification
and $\hinf$-norm estimation have the same worst-case sample complexity.  In the
active setting, the lower bound improves by a factor that is logarithmic in the
filter length. Experimentally, we see that the performance of a simple plugin
$\hinf$-norm estimator is competitive with the proposed active algorithms for
norm estimation in the literature.

Our work raises an interesting question as to whether there 
exists an active algorithm attaining the lower bound, or if instead the
lower bound can be sharpened. In
Section~\ref{sec:active_lb_discussion}, we briefly discussed the technical
hurdles that need to be overcome for both cases.
Beyond resolving the gap between the lower bounds, an interesting question
is how does the sample complexity of both model and norm estimation degrade when
the filter length is unknown. Another direction is to extend the algorithms
and analysis beyond single input single output (SISO) systems.

\section*{Acknowledgements}
We thank Jiantao Jiao for insightful discussions regarding minimax lower bounds for functional estimation and for pointing us to the notes of \cite{wu17}.
We also thank Kevin Jamieson for helpful discussions on multi-armed bandits, and Vaishaal Shankar and Eric Jonas for timely PyWren support.
Finally, we thank Mat{\'{i}}as M{\"{u}}ller for sharing with us the implementation of weighted Thompson Sampling in \cite{muller17bandit}.
This work was generously supported in part by ONR awards N00014-17-1-2191, N00014-17-1-2401, and N00014-18-1-2833, the DARPA Assured Autonomy (FA8750-18-C-0101) and Lagrange (W911NF-16-1-0552) programs, and an Amazon AWS AI Research Award. ST is also supported by a Google PhD fellowship.

\bibliographystyle{abbrv}
\bibliography{IEEEabrv,paper}

\clearpage

\section*{Appendix}

\begin{algorithm}[!h]
\caption{Power Method A \cite{rojas2012analyzing}}
\label{alg:power_method_A}
Input: Normalized $u^{(1)}$.\\
\For{$t = 1$ \upshape to $N$}
{
    Perform the experiment $y^{(t)} = Gu^{(t)} + \eta^{(t)}$.\\
    Create the time-reversed $\tilde y^{(t)}$.\\
    $\mu^{(t)} = \|\tilde y^{(t)}\|_2$.\\
    $u^{(t+1)} = \tilde y^{(t)} / \mu^{(t)}$.\\
    $\hat H_t = \sqrt{\mu^{(t-1)}(u^{(t-1)})^\top \tilde y^{(t)}}.$
}
\Return $\hat H_N$.
\end{algorithm}
\vspace{-0.5cm}
\begin{algorithm}[!h]
\caption{Power Method B \cite{wahlberg2010non}}
\label{alg:power_method_B}
Input: Normalized $u^{(1)}$.\\
\For{$t = 1$ \upshape to $N/2$}
{
    Perform the experiment $y^{(t)} = Gu^{(t)} + \eta^{(t)}$.\\
    Create the time-reversed $\tilde y^{(t)}$.\\
    Perform the experiment $z^{(t)} = G\tilde y^{(t)} + \eta'^{(t)}$.\\
    Create the time-reversed $\tilde z^{(t)}$.\\
    $\hat H_t = \sqrt{|(u^{(t)})^\top\tilde z^{(t)}}|$.\\
    $u^{(t+1)} = \tilde z^{(t)}/\|\tilde z^{(t)}\|_2$.
}
\Return $\hat H_{(N/2)}$.
\end{algorithm}
\vspace{-0.5cm}
\begin{algorithm}[!h]
\caption{Weighted Thompson Sampling (WTS) \cite{muller17bandit}}
\label{alg:WTS}
Input: $M, \lambda, \sigma$, $\rho_k^1 = 1/K\; \forall\; k$, $m^0=0$, $v^0 = \lambda^2 I$.\\
\For{$t = 1$ \upshape to $N$}
{
   Input design: Create the normalized input signal $u^{(t)}$ proportional to the DFT power profile $p^t=\rho^t$.\\
   Perform the experiment $y^{(t)} = Gu^{(t)} + \eta^{(t)}$ and obtain DFT coefficients $X_k^t = Y_k^t / U_k^t$. \\
   Update the posterior for all $k$:\\
   \Indp $m_k^{t+1} = \frac{\lambda^2\sum_{\ell=1}^t p_k^\ell X_k^\ell}{\sigma^2 + \lambda^2\sum_{\ell=1}^t p_k^\ell}$. \\
   $v_k^{t+1} = \lambda^2/(1+\lambda^2/\sigma^2\sum_{\ell=1}^t p_k^\ell)$. \\
   \Indm Update the posterior $\rho^{t+1}$:\\
   \Indp Draw $s^l\sim \mathcal{N_C}(m_k^{t+1},v_k^{t+1}),\quad l=1,\ldots,M$. \\
   $\rho_k^{t+1} = \frac{1}{M}\sum_{l=1}^M \#(\argmax_i\{|s_i^l|\}=k)$. \\
   \Indm $\hat H_t = \max_k \sum_{\ell=1}^t p_k^\ell X_k^\ell / \sum_{\ell=1}^t p_k^\ell$.
}
\Return $\hat H_N$.
\end{algorithm}

\end{document}